\documentclass[11pt,a4paper,twoside]{amsart}

\usepackage{graphics}
\usepackage{color}
\usepackage{epic}
\usepackage{mfpic}

\newtheorem{theo}{\bf Theorem}[section]
\newtheorem{propo}[theo]{\bf Proposition}
\newtheorem{lemma}[theo]{\bf Lemma}
\newtheorem{conj}[theo]{\bf Conjecture}
\newtheorem{defi}[theo]{\bf Definition}

\newcommand{\sh}{{\rm sh}}
\newcommand{\T}{{\mathcal T}}
\newcommand{\SW}{{\rm SW}}
\def\fb{{\rm fb}}
\def\hs{{\rm hs}}
\def\vs{{\rm vs}}

\def\sgn{{\rm\,sgn}}
\def\SYT{{\rm SYT}}

\def\SCT{{\rm SCT}}

\begin{document}

\title{On the sign-imbalance of partition shapes}
\author{Jonas Sj{\"o}strand}
\address{Department of Mathematics, Royal Institute of Technology \\
   SE-100 44 Stockholm, Sweden}
\email{jonass@kth.se}
\date{17 November 2004}

\begin{abstract}
Let the {\em sign} of a standard Young tableau be the sign of the
permutation you get by reading it row by row from left to right,
like a book. A conjecture by Richard Stanley
says that the sum of the signs of all SYTs with $n$ squares
is $2^{\lfloor n/2\rfloor}$. We present a stronger theorem with a
purely combinatorial proof using the Robinson-Schensted correspondence
and a new concept called chess tableaux.

We also prove a sharpening of another conjecture by Stanley concerning
weighted sums of squares of sign-imbalances. The proof is built on a
remarkably simple relation between the sign of a permutation and the
signs of its RS-corresponding tableaux.

\end{abstract}
\maketitle
\section{Introduction}
\noindent
Young tableaux are simple combinatorial objects with complex
properties. They play a central role in the
theory of symmetric functions (see~\cite{fulton}) so
they have been studied a lot,
but the subject is still very much alive.
Recently Richard Stanley came up with a very nice conjecture
on Young tableaux:

{\em Let the {\em sign} of a standard Young tableau be the sign of the
permutation you get by reading it row by row from left to right,
like a book. The sum of the signs of all SYTs with $n$ squares
is $2^{\lfloor n/2\rfloor}$.}

If we take $n=3$ for example, there are four SYTs:

\begin{center}
\setlength{\unitlength}{0.5mm}
\begin{tabular}{c}
$+1$ \\
\begin{picture}(30,40)(0,-35)
\put(0,0){\line(1,0){30}}
\put(0,-10){\line(1,0){30}}
\put(0,0){\line(0,-1){10}}
\put(10,0){\line(0,-1){10}}
\put(20,0){\line(0,-1){10}}
\put(30,0){\line(0,-1){10}}
\put(0,-10){\makebox(10,10){1}}
\put(10,-10){\makebox(10,10){2}}
\put(20,-10){\makebox(10,10){3}}
\end{picture}
\end{tabular}
\begin{tabular}{c}
$+1$ \\
\begin{picture}(30,40)(-5,-35)
\put(0,0){\line(1,0){20}}
\put(0,-10){\line(1,0){20}}
\put(0,-20){\line(1,0){10}}
\put(0,0){\line(0,-1){20}}
\put(10,0){\line(0,-1){20}}
\put(20,0){\line(0,-1){10}}
\put(0,-10){\makebox(10,10){1}}
\put(10,-10){\makebox(10,10){2}}
\put(0,-20){\makebox(10,10){3}}
\end{picture}
\end{tabular}
\begin{tabular}{c}
$-1$ \\
\begin{picture}(30,40)(-5,-35)
\put(0,0){\line(1,0){20}}
\put(0,-10){\line(1,0){20}}
\put(0,-20){\line(1,0){10}}
\put(0,0){\line(0,-1){20}}
\put(10,0){\line(0,-1){20}}
\put(20,0){\line(0,-1){10}}
\put(0,-10){\makebox(10,10){1}}
\put(10,-10){\makebox(10,10){3}}
\put(0,-20){\makebox(10,10){2}}
\end{picture}
\end{tabular}
\begin{tabular}{c}
$+1$ \\
\begin{picture}(30,40)(-10,-35)
\put(0,0){\line(1,0){10}}
\put(0,-10){\line(1,0){10}}
\put(0,-20){\line(1,0){10}}
\put(0,-30){\line(1,0){10}}
\put(0,0){\line(0,-1){30}}
\put(10,0){\line(0,-1){30}}
\put(0,-10){\makebox(10,10){1}}
\put(0,-20){\makebox(10,10){2}}
\put(0,-30){\makebox(10,10){3}}
\end{picture}
\end{tabular}
\end{center}
Their signs sum up to $2=2^{\lfloor 3/2\rfloor}$.

The above conjecture is just a special case of another one
which Stanley gave in~\cite{stanley} (our
conjecture~\ref{conj:stanleyab}(a)).
That conjecture was proved by Lam~\cite{lam} but we will prove
an even stronger theorem (our theorem~\ref{theo:jonas}).
Part (b) of the same conjecture is also proved in
a stronger version (our theorems~\ref{theo:betterb}
and~\ref{theo:special23}).

To settle the conjectures we use two tools: the Robinson-Schensted
correspondence, and a new concept called chess tableaux. Some of
our results in developing these tools have the flavour of an ad hoc lemma,
but proposition~\ref{propo:rscorr}, which is a link between signs of
tableaux and signs of permutations, may be of interest in its own right.

I would like to thank Anders Bj\"orner and Richard Stanley for
introducing me to the ``$2^{\lfloor n/2\rfloor}$-conjecture''.
Many thanks also to an
anonymous referee that has been more than helpful to make this paper
readable.

\section{Preliminaries}
\noindent
An {\bf $n$-shape} $\lambda=(\lambda_1,\lambda_2,\ldots)$ is
a graphical representation (a Ferrers diagram)
of an integer partition of $n=\sum_i\lambda_i$. We write
$\lambda\vdash n$ and we will not distinguish the partition
itself from its shape. The {\em coordinates} of a square
is the pair $(r,c)$ where $r$ and $c$ are the row and column indices.
Example:
\begin{center}
\setlength{\unitlength}{0.5mm}
\begin{picture}(50,60)(0,-55)
\put(-70,-27){$(5,3,2,2,1)\ \ =$}
\put(0,0){\line(1,0){50}}
\put(0,-10){\line(1,0){50}}
\put(0,-20){\line(1,0){30}}
\put(0,-30){\line(1,0){20}}
\put(0,-40){\line(1,0){20}}
\put(0,-50){\line(1,0){10}}
\put(0,0){\line(0,-1){50}}
\put(10,0){\line(0,-1){50}}
\put(20,0){\line(0,-1){40}}
\put(30,0){\line(0,-1){20}}
\put(40,0){\line(0,-1){10}}
\put(50,0){\line(0,-1){10}}
\put(38,-32){$(3,2)$}
\put(35,-30){\vector(-4,1){20}}
\end{picture}
\end{center}

The {\em conjugate} $\lambda'$ of a shape $\lambda$ is the
reflection of $\lambda$ in the main diagonal, i.e.\ exchanging
rows and columns.

A shape $\lambda$ is a {\em subshape} of a shape $\mu$ if
$\lambda_i\leq\mu_i$ for all $i$.
For any subshape $\lambda\subseteq\mu$ the
{\em skew} shape $\mu/\lambda$ is $\mu$ with $\lambda$ deleted.
Example:
\begin{center}
\setlength{\unitlength}{0.5mm}
\begin{picture}(50,60)(0,-50)
\put(-100,-27){$(5,3,2,2,1)/(3,2,2)\ \ =$}
\put(30,0){\line(1,0){20}}
\put(20,-10){\line(1,0){30}}
\put(20,-20){\line(1,0){10}}
\put(0,-30){\line(1,0){20}}
\put(0,-40){\line(1,0){20}}
\put(0,-50){\line(1,0){10}}
\put(0,-30){\line(0,-1){20}}
\put(10,-30){\line(0,-1){20}}
\put(20,-10){\line(0,-1){10}}
\put(20,-30){\line(0,-1){10}}
\put(30,0){\line(0,-1){20}}
\put(40,0){\line(0,-1){10}}
\put(50,0){\line(0,-1){10}}
\end{picture}
\end{center}

A {\bf domino} is a rectangle consisting
of two squares. By $v(\lambda)$ we will
denote the maximal number of disjoint vertical dominoes that fit in the
shape $\lambda$. We let $h(\lambda)=v(\lambda')$.

\begin{figure}
\begin{center}
\resizebox{25mm}{!}{\includegraphics{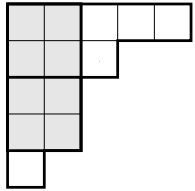}}
\caption{The shaded squares form the fourling body and the white squares
are the strip. Here $d(\lambda)=2$ and $\vs(\lambda)=\hs(\lambda)=1$.}
\label{fig:fourlingbody}
\end{center}
\end{figure}

A {\bf fourling} is a $2\times2$-square. The maximal number
of disjoint fourlings that fit in a shape $\lambda$ is denoted by
$d(\lambda)$. A {\bf fourling shape} is a (possibly empty)
shape consisting of fourlings.
The {\bf fourling body} $\fb(\lambda)$ of a shape $\lambda$
is its largest fourling subshape.
The remaining squares form the
{\bf strip} of the shape. By $\vs(\lambda)$ we will denote
the maximal number of disjoint vertical dominoes that fit in the strip of
$\lambda$. We let $\hs(\lambda)=\vs(\lambda')$.
See figure~\ref{fig:fourlingbody}.

A {\bf tableau} on an $n$-shape $\lambda$ is a labelling of the squares
of $\lambda$ with $n$ different integers such that every integer is greater
than its neighbours above and to the left.
A {\bf standard Young tableau (SYT)} on an $n$-shape
is a tableau with the numbers $[n]=\{1,2,\ldots,n\}$.
We let $\SYT(\lambda)$ denote the set of SYTs on the shape $\lambda$.
Here is an example:

\begin{center}
\setlength{\unitlength}{0.5mm}
\begin{picture}(50,60)(0,-55)
\put(0,0){\line(1,0){50}}
\put(0,-10){\line(1,0){50}}
\put(0,-20){\line(1,0){30}}
\put(0,-30){\line(1,0){20}}
\put(0,-40){\line(1,0){20}}
\put(0,-50){\line(1,0){10}}
\put(0,0){\line(0,-1){50}}
\put(10,0){\line(0,-1){50}}
\put(20,0){\line(0,-1){40}}
\put(30,0){\line(0,-1){20}}
\put(40,0){\line(0,-1){10}}
\put(50,0){\line(0,-1){10}}
\put(0,-10){\makebox(10,10){1}}
\put(10,-10){\makebox(10,10){4}}
\put(20,-10){\makebox(10,10){6}}
\put(30,-10){\makebox(10,10){7}}
\put(40,-10){\makebox(10,10){10}}
\put(0,-20){\makebox(10,10){2}}
\put(10,-20){\makebox(10,10){5}}
\put(20,-20){\makebox(10,10){9}}
\put(0,-30){\makebox(10,10){3}}
\put(10,-30){\makebox(10,10){11}}
\put(0,-40){\makebox(10,10){8}}
\put(10,-40){\makebox(10,10){13}}
\put(0,-50){\makebox(10,10){12}}
\end{picture}
\end{center}
The shape of a tableau $T$ is denoted by $\sh(T)$.

By a {\bf $k$-word} we will mean a sequence of $k$ integers, all different.
A {\bf sorted word} is a strictly increasing sequence of integers.
The {\em sign} of a word $w=w_1w_2\cdots w_k$
is $(-1)^{|\{(i,j)\,:\,i<j,\,w_i>w_j\}|}$, so it is
$+1$ for an even number of inversions, $-1$ otherwise.

The {\bf sign} $\sgn(T)$ of a tableau $T$ is
the sign of the word you get
by reading the integers row by row, from left to right and from
top to bottom, like a book.
Our example tableau has 18 inversions, so $\sgn(T)=+1$.
The {\bf sign-imbalance} $I_\lambda$ of a shape
$\lambda$ is the sum of the signs of all SYTs on that shape.
\begin{defi}
$$I_\lambda=\sum_{T\in\SYT(\lambda)}\sgn(T).$$
\end{defi}

\section{Stanley's conjecture and our results}
\noindent
Richard Stanley gave the following conjecture in \cite{stanley}.
\begin{conj}\label{conj:stanleyab}
\mbox{\\}
\begin{enumerate}
\item[(a)] For every $n\geq0$
$$\sum_{\lambda\vdash n}q^{v(\lambda)}t^{d(\lambda)}x^{h(\lambda)}
I_\lambda=(q+x)^{\lfloor n/2\rfloor}.$$
\item[(b)] If $n\not\equiv1\pmod4$
$$\sum_{\lambda\vdash n}(-1)^{v(\lambda)}t^{d(\lambda)}I_\lambda^2=0.$$
\end{enumerate}
\end{conj}

\noindent
The special case $t=0$ of (a) goes like this:
\begin{propo}\label{propo:hooks}
For all $n\geq0$ we have
$$\sum_{\lambda=(n-i,1^i)}q^{v(\lambda)}x^{h(\lambda)}I_\lambda
=(q+x)^{\lfloor n/2\rfloor},$$
where $\lambda$ ranges over all hooks $(n-i,1^i)$, $0\leq i\leq n-1$.
\end{propo}
\noindent
It tells us that the right hand side
$(q+x)^{\lfloor n/2\rfloor}$ comes from the {\em hooks}, i.e
the fourling-free shapes, and was proved twice by Stanley in
\cite[prop.\ 3.4]{stanley}. We give a third proof in section~\ref{sec:proof}.

The rest of (a) says that,
for fixed $d\geq1$, $h$ and $v$,
the sum of the sign-imbalances of all $n$-shapes $\lambda$
with $v(\lambda)=v$, $h(\lambda)=h$ and $d(\lambda)=d$ vanishes.

Part (a) of the conjecture has been proved by Lam \cite{lam}.
We will prove a stronger version of part (a)
which lets us fix not only the number of fourlings
but the whole fourling shape:
\begin{theo}\label{theo:jonas}
Given a non-empty fourling shape $D$
and nonnegative integers $h$, $v$ and $s$,
$$\sum I_\lambda=0$$
where the sum is taken over all shapes $\lambda$ with fourling body
$D$, $s$ squares in the strip, $\hs(\lambda)=h$, and $\vs(\lambda)=v$.
\end{theo}
\noindent
The proof will be found in section~\ref{sec:mainproof} and is purely
combinatorial. Figure~\ref{fig:jonasex}
\begin{figure}
\begin{center}
\setlength{\unitlength}{0.3mm}
\begin{tabular}{c}
5 \\
\begin{picture}(90,60)(-10,-55)
\thicklines
\put(0,0){\line(1,0){40}}
\put(0,-20){\line(1,0){40}}
\put(0,0){\line(0,-1){20}}
\put(40,0){\line(0,-1){20}}
\thinlines
\put(0,0){\line(1,0){70}}
\put(0,-10){\line(1,0){70}}
\put(0,-20){\line(1,0){50}}
\put(0,0){\line(0,-1){20}}
\put(10,0){\line(0,-1){20}}
\put(20,0){\line(0,-1){20}}
\put(30,0){\line(0,-1){20}}
\put(40,0){\line(0,-1){20}}
\put(50,0){\line(0,-1){20}}
\put(60,0){\line(0,-1){10}}
\put(70,0){\line(0,-1){10}}
\end{picture}
\end{tabular}
\begin{tabular}{c}
5 \\
\begin{picture}(90,60)(-15,-55)
\thicklines
\put(0,0){\line(1,0){40}}
\put(0,-20){\line(1,0){40}}
\put(0,0){\line(0,-1){20}}
\put(40,0){\line(0,-1){20}}
\thinlines
\put(0,0){\line(1,0){60}}
\put(0,-10){\line(1,0){60}}
\put(0,-20){\line(1,0){50}}
\put(0,-30){\line(1,0){10}}
\put(0,0){\line(0,-1){30}}
\put(10,0){\line(0,-1){30}}
\put(20,0){\line(0,-1){20}}
\put(30,0){\line(0,-1){20}}
\put(40,0){\line(0,-1){20}}
\put(50,0){\line(0,-1){20}}
\put(60,0){\line(0,-1){10}}
\end{picture}
\end{tabular}
\begin{tabular}{c}
2 \\
\begin{picture}(90,60)(-15,-55)
\thicklines
\put(0,0){\line(1,0){40}}
\put(0,-20){\line(1,0){40}}
\put(0,0){\line(0,-1){20}}
\put(40,0){\line(0,-1){20}}
\thinlines
\put(0,0){\line(1,0){60}}
\put(0,-10){\line(1,0){60}}
\put(0,-20){\line(1,0){40}}
\put(0,-30){\line(1,0){10}}
\put(0,-40){\line(1,0){10}}
\put(0,0){\line(0,-1){40}}
\put(10,0){\line(0,-1){40}}
\put(20,0){\line(0,-1){20}}
\put(30,0){\line(0,-1){20}}
\put(40,0){\line(0,-1){20}}
\put(50,0){\line(0,-1){10}}
\put(60,0){\line(0,-1){10}}
\end{picture}
\end{tabular}
\begin{tabular}{c}
2 \\
\begin{picture}(90,60)(-20,-55)
\thicklines
\put(0,0){\line(1,0){40}}
\put(0,-20){\line(1,0){40}}
\put(0,0){\line(0,-1){20}}
\put(40,0){\line(0,-1){20}}
\thinlines
\put(0,0){\line(1,0){50}}
\put(0,-10){\line(1,0){50}}
\put(0,-20){\line(1,0){50}}
\put(0,-30){\line(1,0){20}}
\put(0,0){\line(0,-1){30}}
\put(10,0){\line(0,-1){30}}
\put(20,0){\line(0,-1){30}}
\put(30,0){\line(0,-1){20}}
\put(40,0){\line(0,-1){20}}
\put(50,0){\line(0,-1){20}}
\end{picture}
\end{tabular}
\begin{tabular}{c}
$-5$ \\
\begin{picture}(90,60)(-20,-55)
\thicklines
\put(0,0){\line(1,0){40}}
\put(0,-20){\line(1,0){40}}
\put(0,0){\line(0,-1){20}}
\put(40,0){\line(0,-1){20}}
\thinlines
\put(0,0){\line(1,0){50}}
\put(0,-10){\line(1,0){50}}
\put(0,-20){\line(1,0){40}}
\put(0,-30){\line(1,0){20}}
\put(0,-40){\line(1,0){10}}
\put(0,0){\line(0,-1){40}}
\put(10,0){\line(0,-1){40}}
\put(20,0){\line(0,-1){30}}
\put(30,0){\line(0,-1){20}}
\put(40,0){\line(0,-1){20}}
\put(50,0){\line(0,-1){10}}
\end{picture}
\end{tabular}
\begin{tabular}{c}
$-2$ \\
\begin{picture}(90,60)(-25,-55)
\thicklines
\put(0,0){\line(1,0){40}}
\put(0,-20){\line(1,0){40}}
\put(0,0){\line(0,-1){20}}
\put(40,0){\line(0,-1){20}}
\thinlines
\put(0,0){\line(1,0){40}}
\put(0,-10){\line(1,0){40}}
\put(0,-20){\line(1,0){40}}
\put(0,-30){\line(1,0){30}}
\put(0,-40){\line(1,0){10}}
\put(0,0){\line(0,-1){40}}
\put(10,0){\line(0,-1){40}}
\put(20,0){\line(0,-1){30}}
\put(30,0){\line(0,-1){30}}
\put(40,0){\line(0,-1){20}}
\end{picture}
\end{tabular}
\begin{tabular}{c}
$-7$ \\
\begin{picture}(90,60)(-25,-55)
\thicklines
\put(0,0){\line(1,0){40}}
\put(0,-20){\line(1,0){40}}
\put(0,0){\line(0,-1){20}}
\put(40,0){\line(0,-1){20}}
\thinlines
\put(0,0){\line(1,0){40}}
\put(0,-10){\line(1,0){40}}
\put(0,-20){\line(1,0){40}}
\put(0,-30){\line(1,0){20}}
\put(0,-40){\line(1,0){10}}
\put(0,-50){\line(1,0){10}}
\put(0,0){\line(0,-1){50}}
\put(10,0){\line(0,-1){50}}
\put(20,0){\line(0,-1){30}}
\put(30,0){\line(0,-1){20}}
\put(40,0){\line(0,-1){20}}
\end{picture}
\end{tabular}
\setlength{\unitlength}{0.2mm}
\caption[trams]{The imbalances of the 12-shapes $\lambda$ with fourling
body
\begin{picture}(40,20)(0,-15)
\put(0,0){\line(1,0){40}}
\put(0,-10){\line(1,0){40}}
\put(0,-20){\line(1,0){40}}
\put(0,0){\line(0,-1){20}}
\put(10,0){\line(0,-1){20}}
\put(20,0){\line(0,-1){20}}
\put(30,0){\line(0,-1){20}}
\put(40,0){\line(0,-1){20}}
\end{picture}
and $\vs(\lambda)=\hs(\lambda)=1$. You can check that their sum vanishes.}
\label{fig:jonasex}
\end{center}
\end{figure}
shows an example.

In the same spirit, we have the following theorem which
is a sharpening of (b) when $n$ is even.

\begin{theo}\label{theo:betterb}
Given a fourling shape $D$ and an even integer $n\geq0$,
$$\sum(-1)^{v(\lambda)}I_\lambda^2=0$$
where the sum is taken over all $n$-shapes $\lambda$ with
$\fb(\lambda)=D$.
\end{theo}
\noindent
We will prove it in section~\ref{sec:rs}.

The next theorem, which we prove in section~\ref{sec:chess},
covers the rest of (b).
\begin{theo}\label{theo:special23}
If $n\equiv2$ or $n\equiv3\pmod4$
$$\sum_{\lambda\vdash n}(-1)^{v(\lambda)}F(\lambda)=0$$
for any function
$F:\{\mbox{$n$-shapes}\}\rightarrow\mathbb{C}$
such that $F(\lambda)=F(\lambda')$
and $I_\lambda=0\Rightarrow F(\lambda)=0$
for all $n$-shapes $\lambda$.
\end{theo}
\noindent
Choosing $F(\lambda)=t^{d(\lambda)}I_\lambda^2$ proves (b) for
$n\equiv2$ and $n\equiv3\pmod4$ since $|I_\lambda|=|I_\lambda'|$
(see e.g.\ Stanley \cite{stanley} or our proposition~\ref{propo:transpose}).
Thus we have proved all parts
of Stanley's conjecture.

Finally, the special case $t=1$ of (b) will be proved
also without the assumption $n\not\equiv1\pmod4$:
\begin{theo}\label{theo:special}
For all $n\geq2$
$$\sum_{\lambda\vdash n}(-1)^{v(\lambda)}I_\lambda^2=0.$$
\end{theo}
\noindent
This was proved independently of us by Reifegerste \cite[theorem 5.1]{astrid}.
Stanley proved it for even $n$ \cite[theorem 3.2(b)]{stanley}.

The rest of this paper is composed as follows.
In section~\ref{sec:chess} we introduce the concept of a chess tableau
and prove theorem~\ref{theo:special23}. In section~\ref{sec:rs}
we show how the signs of tableaux and permutations are
related by the Robinson-Schensted correspondence. The most important
result is proposition~\ref{propo:rscorr} which we use
to prove theorem~\ref{theo:special} and~\ref{theo:betterb}.
Finally, in section~\ref{sec:proof} we prove theorem~\ref{theo:jonas}
using chess tableaux and the RS-correspondence.

\section{Chess tableaux and theorem~\ref{theo:special23}\label{sec:chess}}
\noindent
When working on sums of tableau signs one is naturally led
to use domino tableaux (see~\cite{stanley} and~\cite{white}).
In this paper we choose a similar approach which turns out to be
more successful in settling the conjectures.

A {\em chess colouring} of a shape is a colouring of the squares such that
a square $(r,c)$ is black if $r+c$ is even and white if $r+c$ is odd.
From now on we will frequently refer to
white and black squares of a shape, implicitly meaning the chess colouring.
A {\bf chess tableau} is a SYT with odd integers in
black squares and even in white.

\begin{lemma}\label{lemma:chess}
Given a shape $\lambda$, $\sum_{T\in\SCT(\lambda)}\sgn(T)=I_\lambda$,
where $\SCT(\lambda)$ is the set of chess tableaux on $\lambda$.
\end{lemma}
\begin{proof}
There is a sign-alternating involution on the
non-chess SYTs: Given a non-chess SYT there are at least two
consecutive integers of the same colour. Choose the least such pair
and switch the integers. This is allowed unless they are
horizontal or vertical neighbours, which they are not since
neighbours have different colours.
\end{proof}

\begin{propo}\label{propo:colourbalance}
If $\lambda$ is a shape with $s$ strip squares,
$I_\lambda\neq0$ only if it has equally many white and black squares
or one more black square. This implies that
$\hs(\lambda)+\vs(\lambda)=\lfloor s/2\rfloor$.
\end{propo}
\begin{proof}
Let $B$ and $W$ be the number of
black respectively white squares in the strip of $\lambda$.
By lemma~\ref{lemma:chess} we must have $B=W$ or
$B=W+1$ if $I_\lambda\neq0$ (otherwise there are no chess tableaux).
Every white strip square belongs
to a certain strip domino, namely the one with the black square
above or to the left,
so $W=\hs(\lambda)+\vs(\lambda)$, see
figure~\ref{fig:chess}.
\begin{figure}
\begin{center}
\resizebox{25mm}{!}{\includegraphics{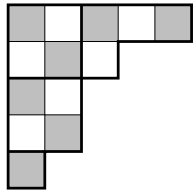}}
\caption{The white strip squares count the strip dominoes, $\vs(\lambda)+\hs(\lambda)=2$.}
\label{fig:chess}
\end{center}
\end{figure}
Thus, for a $\lambda$
with $I_\lambda\neq0$ we have $\hs(\lambda)+\vs(\lambda)=\lfloor s/2\rfloor$.
\end{proof}

\begin{proof}[\bf Proof of theorem~\ref{theo:special23}]
We show that
if $\lambda$ is an $n$-shape with $n\equiv2$ or $n\equiv3\pmod4$,
either $I_\lambda=0$ or $v(\lambda)\not\equiv h(\lambda)\pmod2$.
This implies that the non-vanishing terms $(-1)^{v(\lambda)}F(\lambda)$
come in cancelling pairs
$(-1)^{v(\lambda)}F(\lambda)+(-1)^{v(\lambda')}F(\lambda')$.

Suppose $I_\lambda\neq0$ and
let $s$ be the number of strip squares in $\lambda$.
Since the fourling body consists of fourlings
we have $s\equiv2$ or $s\equiv3\pmod4$.
By proposition~\ref{propo:colourbalance} we can assume that
$\hs(\lambda)+\vs(\lambda)=\lfloor s/2\rfloor$ which is odd.
The fourling body has equally many horizontal and vertical dominoes
so $v(\lambda)\not\equiv h(\lambda)\pmod2$.
\end{proof}

\section{Robinson-Schensted correspondence
and theorems~\ref{theo:special} and~\ref{theo:betterb}\label{sec:rs}}
\noindent
Given a tableau $T$ and a number $a$ different from all
numbers in $T$, by
{\bf (row) insertion} of $a$ into $T$
we mean the usual Robinson-Schensted insertion
(see for example \cite[p.\ 316]{enum2})
resulting in a tableau $(T\leftarrow a)$ with one more square $x$
than $T$. By {\bf (row) extraction} of $x$ we mean
the reverse process resulting in $T$ and $a$.
Insertion of a word into a tableau means
insertion of the integers in the word one by one from left to right.

We will use the following lemma later on.
\begin{lemma}\label{lemma:noncrossing}
Given a tableau $T$ and integers $a\neq b$ different from all
entries in $T$, the square
$\sh(T\leftarrow ab)/\sh(T\leftarrow a)$ appears in a column
somewhere to the right of $\sh(T\leftarrow a)/\sh(T)$
if and only if $a<b$.
\end{lemma}

\begin{proof}
Suppose that $a<b$.
We can insert the two numbers in parallel row by row. If $a$
is greater than every number in the first row, the squares
$x=\sh(T\leftarrow a)/\sh(T)$ and
$y=\sh(T\leftarrow ab)/\sh(T\leftarrow a)$ will be
placed rightmost in that row with $y$ to the right of $x$.
If $a$ pops a number $a_2$ in the first row, $b$ will either
terminate leaving $y$ rightmost in the first row or pop a number
$b_2>a_2$. The if part of the lemma follows by induction.
The converse is proved similarly.
\end{proof}

The next lemma tells us what insertion does to the sign of
the tableau.
\begin{lemma}\label{lemma:insertone}
If $T$ is a tableau and $a$ is a number different from all entries in $T$,
$$\sgn(T\leftarrow a)=(-1)^{l+w+u}\sgn(T),$$
where $l$ is the number of entries
in $T$ less than $a$, $w$ is 0 if $\sh(T\leftarrow a)/\sh(T)$
is black and 1 if it is white,
and $u$ is the number of squares in rows above $\sh(T\leftarrow a)/\sh(T)$.
\end{lemma}
\begin{figure}
\begin{center}
\resizebox{50mm}{!}{\includegraphics{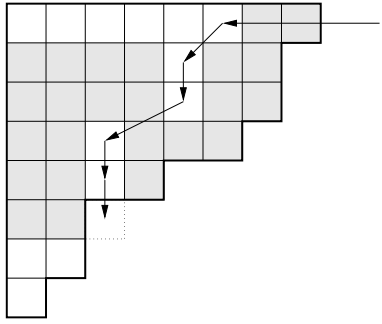}}
\caption{Insertion of a number. The shaded squares are
counted by $\sum_{i=2}^k(\lambda_{i-1}-c_{i-1}+c_i-1)$ in the proof.}
\label{fig:rsinsert}
\end{center}
\end{figure}
\begin{proof}
Let $\lambda=\sh(T)$ and look at figure~\ref{fig:rsinsert}.
During the insertion $a_1=a$ pops a number $a_2$ at $(1,c_1)$ which pops
a number $a_3$ at $(2,c_2)$ and so on. Finally the number $a_k$ fills
a new square $(k,c_k)=\sh(T\leftarrow a)/\sh(T)$.
For $2\leq i\leq k$, the move of $a_i$ multiplies
the sign of the tableau by $(-1)^{\lambda_{i-1}-c_{i-1}+c_i-1}$.
Summation yields
$$\sum_{i=2}^k(\lambda_{i-1}-c_{i-1}+c_i-1)
=c_k-c_1+\sum_{i=1}^{k-1}(\lambda_i-1)=u-k+1+c_k-c_1.$$
The placing of $a=a_1$ in the first row multiplies the sign of the tableau by
$(-1)^{l-c_1+1}$, so the total factor is
$(-1)^{u-k+1+c_k-c_1+l-c_1+1}=(-1)^{u+l+c_k+k}=(-1)^{u+l+w}$.
\end{proof}
Now the following
natural question arises: How is the sign property transferred
by the RS-correspondence? The answer is quite beautiful:
\begin{propo}\label{propo:rscorr}
In the RS-correspondence $\pi\leftrightarrow(P,Q)$ we have
$$\sgn(\pi)=(-1)^{v(\lambda)}\sgn(P)\sgn(Q)$$
where $\lambda$ is the shape of $P$ and $Q$.
\end{propo}
\begin{proof}
Suppose we have inserted the first $k$ numbers in $\pi$ yielding
tableaux $P^k$ and $Q^k$ on the shape $\lambda^k$, and
$\sgn(\pi_1\cdots\pi_k)=(-1)^{v(\lambda^k)}\sgn(P^k)\sgn(Q^k)$.
This is certainly true for $k=0$. Now we argue by induction
over $k$.
We insert the next number $\pi_{k+1}$ and look what happens
according to lemma~\ref{lemma:insertone}.
We get
$\sgn(P^{k+1})=(-1)^{l+w+u}\sgn(P^k)$, and if $\lambda^{k+1}/\lambda^k$
has coordinates $(r,c)$ we get
$\sgn(Q^{k+1})=(-1)^{k-u-c+1}\sgn(Q^k)=(-1)^{k-u-w+r+1}\sgn(Q^k)$
since $w$ is congruent to $r+c$ modulo 2.
Whether a new vertical domino will fit in $\lambda^{k+1}$ is
only dependent on $r$, so
$(-1)^{v(\lambda^{k+1})}=(-1)^{r+1}(-1)^{v(\lambda^k)}$.
Finally, $\sgn(\pi_1\cdots\pi_{k+1})=(-1)^{k-l}\sgn(\pi_1\cdots\pi_k)$.

Putting it all together yields at last
$$\sgn(\pi_1\cdots\pi_{k+1})=(-1)^{k-l}\sgn(\pi_1\cdots\pi_k)
=(-1)^{k-l}(-1)^{v(\lambda^k)}\sgn(P^k)\sgn(Q^k)=$$
$$=(-1)^{r+1}(-1)^{v(\lambda^k)}
(-1)^{l+w+u}\sgn(P^k)(-1)^{k-u-w+r+1}\sgn(Q^k)=$$
$$=(-1)^{v(\lambda^{k+1})}\sgn(P^{k+1})\sgn(Q^{k+1}).$$
\end{proof}
\noindent
The above result was also found by Reifegerste
\cite[theorem 4.3]{astrid} independently of us.
\vspace{2mm}
\paragraph{\bf Remark.}
If we specialise to the RS-bijection $\pi\leftrightarrow(P,P)$
between involutions $\pi\in S_n$ and $n$-SYTs $P$,
proposition~\ref{propo:rscorr} gives
that $\sgn(\pi)=(-1)^{v(\sh(P))}$. This is also a simple consequence
of a theorem by Sch\"utzenberger \cite[page 127]{schutzenberger}
(see also~\cite[exercise~7.28~a]{enum2})
stating that the number of
fix points in $\pi$ equals the number of columns of $P$ of
odd length.
\vspace{2mm}

As a simple consequence of proposition~\ref{propo:rscorr}
we get theorem~\ref{theo:special}.
\begin{proof}[\bf Proof of theorem~\ref{theo:special}]
By proposition~\ref{propo:rscorr} we have
$$\sum_{\lambda\vdash n}(-1)^{v(\lambda)}I_\lambda^2
=\sum_{\lambda\vdash n}(-1)^{v(\lambda)}
\left(\sum_{P\in\SYT(\lambda)}\sgn(P)\right)^2=$$
$$=\sum_{\lambda\vdash n}\sum_{P,Q\in\SYT(\lambda)}
(-1)^{v(\lambda)}\sgn(P)\sgn(Q)
=\sum_{\pi\in S_n}\sgn(\pi)=0.$$
\end{proof}
To prove theorem~\ref{theo:betterb} we will need the following much stronger
theorem which is proved in a manner similar to what we did above.
\begin{theo}\label{theo:bettererb}
Given a set $B$ of black squares and an even integer
$n\geq0$,
$$\sum(-1)^{v(\lambda)}I_\lambda^2=0$$
where the sum is taken over all $n$-shapes $\lambda$
whose black squares are exactly the ones in $B$.
\end{theo}

\begin{proof}
Let $A$ be the set of shapes
whose black squares are exactly the ones in $B$.
For an $n$-SYT $Q$, let $Q\setminus n$ denote the $(n-1)$-SYT we get by
deleting the number $n$ from $Q$.
If $Q$ is a chess tableau,
$\sh(Q)\in A\Leftrightarrow\sh(Q\setminus n)\in A$
since $\sh(Q)$ and $\sh(Q\setminus n)$ contain exactly the
same set of black squares (remember that $n$ is even).
Then, by lemma~\ref{lemma:chess},
$$\sum_{
\mbox{\scriptsize$\begin{array}{c}
\lambda\vdash n \\
\lambda\in A\end{array}$}
}
(-1)^{v(\lambda)}I_\lambda^2
=\sum_{\lambda\vdash n}(-1)^{v(\lambda)}I_\lambda
\sum_{\mbox{\scriptsize$\begin{array}{c}
Q\in\SCT(\lambda)\\
\sh(Q\setminus n)\in A\end{array}$}}\sgn(Q).
$$
Now we take any $n$-shape $\lambda$ and compute its contribution to
the sum.
If $\lambda$ does not have equally many white and black squares,
$I_\lambda=0$ by proposition~\ref{propo:colourbalance} and the
contribution is zero. If $\lambda$ has equally many white and black squares,
then, for $Q\in\SYT(\lambda)$, $Q$ is a chess tableau if and only if
$Q\setminus n$ is a chess tableau. Thus, we can write our expression
in a slightly different way:
$$
\sum_{\lambda\vdash n}(-1)^{v(\lambda)}I_\lambda
\sum_{\mbox{\scriptsize$\begin{array}{c}
Q\in\SYT(\lambda)\\
Q\setminus n\mbox{\ is a chess tableau}\\
\sh(Q\setminus n)\in A\end{array}$}}\sgn(Q)
$$
By proposition~\ref{propo:rscorr} this equals
$$
\sum_{\lambda\vdash n}
\sum_{\mbox{\scriptsize$\begin{array}{c}
P,Q\in\SYT(\lambda)\\
Q\setminus n\mbox{\ is a chess tableau}\\
\sh(Q\setminus n)\in A\end{array}$}}(-1)^{v(\lambda)}\sgn(P)\sgn(Q)
=\sum_{\pi\in S}\sgn(\pi)
$$
where $S\subseteq S_n$ is the set of permutations corresponding
to $n$-tableaux $P$ and $Q$ such that $Q\setminus n$
is a chess tableau whose shape is in $A$. (Note that we do not
require that $Q$ is a chess tableau.)

For an $n$-permutation $\pi$, let $\pi'$ be the
$(n-1)$-permutation defined by
$$\pi'_i=\left\{\begin{array}{lr}
\pi_i & \mbox{if $\pi_i<\pi_n$} \\
\pi_i-1 & \mbox{if $\pi_i>\pi_n$}
\end{array}\right..$$
We can consider the set $S_n$ of $n$-permutations as a disjoint union
$S_n=\bigcup_{\rho\in S_{n-1}}S_n^{\rho}$, where
$S_n^{\rho}=\{\pi\in S_n:\pi'=\rho\}$.
In the RS-correspondence $\pi\rightarrow(P,Q)$ the locations
of the first $n-1$ numbers in $Q$ are only dependent on $\pi'$.
Thus we can write $S$ as a disjoint union
$S=\bigcup_{\rho\in S'}S_n^{\rho}$ where $S'$ is the set of
$(n-1)$-permutations corresponding to a chess $Q$-tableau whose
shape is in $A$.
But $\sum_{\pi\in S_n^{\rho}}\sgn(\pi)=0$ since we can choose the
last element $\pi_n$ in an even number of ways.
\end{proof}

Finally we show that theorem~\ref{theo:betterb}
is a simple consequence of the above theorem.
\begin{proof}[\bf Proof of theorem~\ref{theo:betterb}]
Note that it is impossible to change
the fourling body of a shape by adding or removing
only white squares.

Let $B_\lambda$ denote the set of black squares in a
shape $\lambda$ and let
$\mathcal B=\{B_\lambda\,:\,\lambda\vdash n,\,\fb(\lambda)=D\}$.
Then
$$\sum_{\mbox{\scriptsize
$\begin{array}{c}\lambda\vdash n\\ \fb(\lambda)=D
\end{array}$}}(-1)^{v(\lambda)}I_\lambda^2
=\sum_{B\in\mathcal B}\sum_{\mbox{\scriptsize
$\begin{array}{c}\lambda\vdash n\\ B_\lambda=B
\end{array}$}}(-1)^{v(\lambda)}I_\lambda^2=0$$
by theorem~\ref{theo:bettererb}.
\end{proof}

\section{The proofs of proposition~\ref{propo:hooks}
and theorem~\ref{theo:jonas}\label{sec:proof}}
\label{sec:mainproof}
\noindent
First some definitions:
\begin{defi}\label{defi:sigmatau}
Given an $n$-shape $\lambda$ and an integer $k\geq0$,
let $\T_{\lambda,k}$
be the set of tableaux on $\lambda$ with
numbers in $[n+k]$.

Given $T\in\T_{\lambda,k}$, let {\em the
complementary $k$-word $w_{T,k}$ of $T$} be the sorted $k$-word of
the elements of $[n+k]$ not in $T$.

Let $\SW_{i,j}$ denote the set of sorted $j$-words with
letters in $[i]$.

Given a $k$-word $w$, let $\sigma(w)=(-1)^L$, where $L=\sum_{i=1}^k(w_i-1)$.

Given a skew shape $\mu/\lambda$,
let $\tau(\mu/\lambda)=(-1)^{W+U}$,
where $W$ is the number of white squares in
$\mu/\lambda$ and $U$ is the number of square pairs
$(x,y)\in\lambda\times\mu/\lambda$ with $x$ in a row 
somewhere above $y$.
\end{defi}

\begin{figure}
\fbox{\begin{minipage}{5.5in}
Let $\lambda=(5,2,2,1)$ and $k=3$.
\setlength{\unitlength}{0.5mm}
If we take, for example,
\begin{center}
\begin{picture}(90,40)(-20,-40)
\put(-20,-20){$T=$}
\put(-0.1730,-39.8270){\resizebox{25.3459mm}{!}{\includegraphics{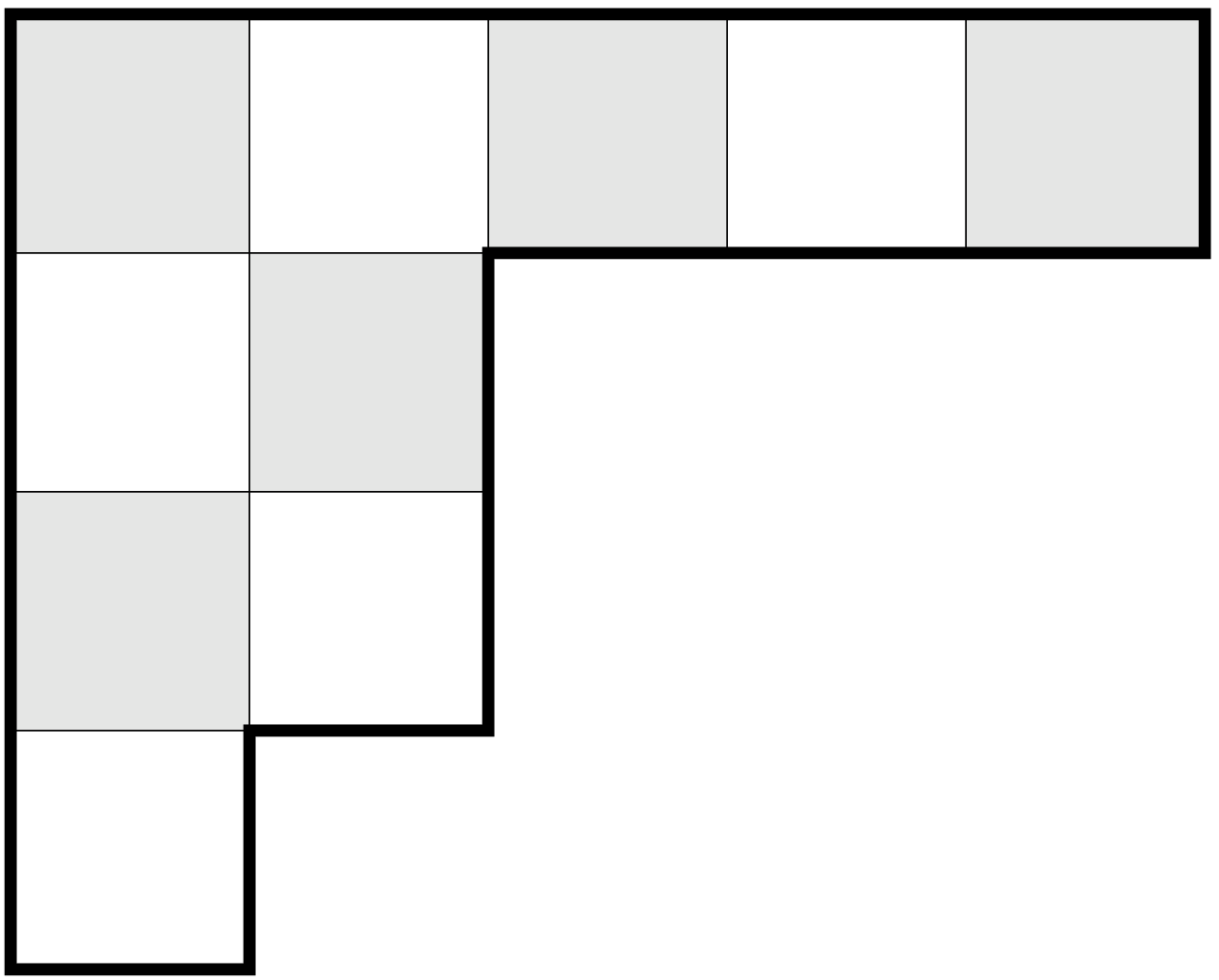}}}
\put(0,-10){\makebox(10,10){2}}
\put(10,-10){\makebox(10,10){4}}
\put(20,-10){\makebox(10,10){6}}
\put(30,-10){\makebox(10,10){9}}
\put(40,-10){\makebox(10,10){13}}
\put(0,-20){\makebox(10,10){3}}
\put(10,-20){\makebox(10,10){5}}
\put(0,-30){\makebox(10,10){8}}
\put(10,-30){\makebox(10,10){11}}
\put(0,-40){\makebox(10,10){12}}
\end{picture}
\end{center}
then $w_{T,3}=1\ 7\ 10$ and insertion yields
\begin{center}
\begin{picture}(160,50)(-55,-50)
\put(-0.1730,-49.8270)
{\resizebox{25.3459mm}{!}{\includegraphics{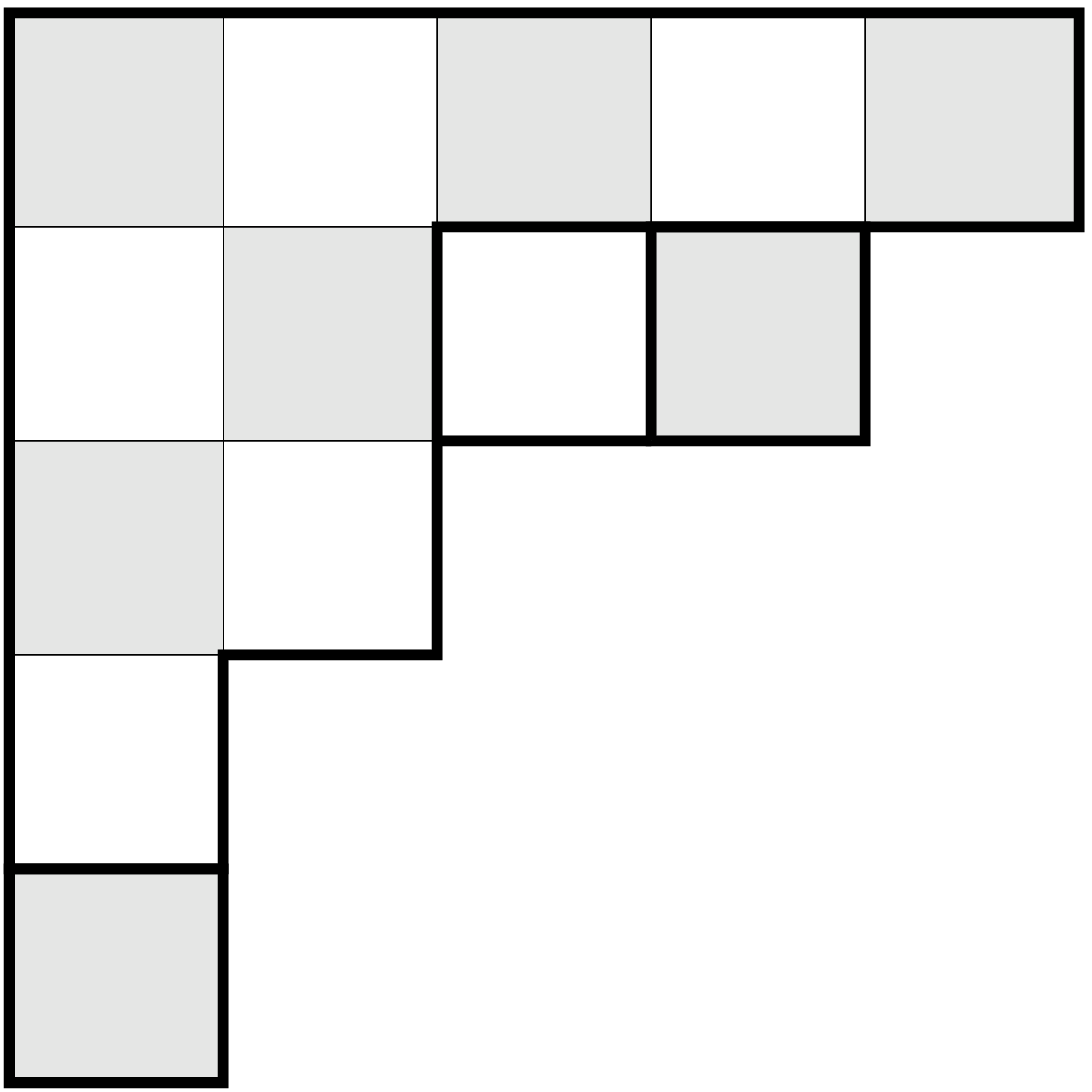}}}
\put(-55,-25){$(T\leftarrow w_{T,3})=$}
\put(0,-10){\makebox(10,10){1}}
\put(10,-10){\makebox(10,10){4}}
\put(20,-10){\makebox(10,10){6}}
\put(30,-10){\makebox(10,10){7}}
\put(40,-10){\makebox(10,10){10}}
\put(0,-20){\makebox(10,10){2}}
\put(10,-20){\makebox(10,10){5}}
\put(20,-20){\makebox(10,10){9}}
\put(30,-20){\makebox(10,10){13}}
\put(0,-30){\makebox(10,10){3}}
\put(10,-30){\makebox(10,10){11}}
\put(0,-40){\makebox(10,10){8}}
\put(0,-50){\makebox(10,10){12}}
\end{picture}
\end{center}
We get $L=(1-1)+(7-1)+(10-1)=15$, so
$\sigma(1\ 7\ 10)=(-1)^L=-1$. Among the three extra squares
only one is white, so $W=1$. The number of original squares
in rows above the extra squares is 10, 5 and 5, so $U=20$, and
$\tau(\sh(T\leftarrow w_{T,3})/\lambda)=(-1)^{W+U}=(-1)^{1+20}=-1$.
The lemma says that $\sgn(T\leftarrow w_{T,3})=\sigma\tau\sgn(T)$.
We check that $T$ has 11 inversions and $(T\leftarrow w_{T,3})$
has 21, so it seems alright.
\end{minipage}}
\caption[trams]{\label{fig:wordinsertion}
Example of lemma~\ref{lemma:wordinsertion}.}
\end{figure}

\begin{lemma}\label{lemma:wordinsertion}
Let $\lambda$ be an $n$-shape.
Insertion of $w_{T,k}$ into $T$ gives a bijection between
$\T_{\lambda,k}$ and the set of SYTs on $(n+k)$-shapes
$\mu\supseteq\lambda$ with $v(\mu/\lambda)=0$. We have
\begin{equation}\label{eq:wordinsertion}
\sgn(T\leftarrow w_{T,k})
=\sigma(w_{T,k})\tau(\sh(T\leftarrow w_{T,k})/\lambda)\sgn(T).
\end{equation}
Figure~\ref{fig:wordinsertion} shows an example.
\end{lemma}
\begin{proof}
Let $T\in\T_{\lambda,k}$ and let
$\mu=\sh(T\leftarrow w_{T,k})$. By lemma~\ref{lemma:noncrossing}
the extra squares $\mu/\lambda$ will appear from left to right,
without any vertical dominoes.
The inverse of the insertion is extraction
of the squares $\mu/\lambda$ from right to left. Clearly it
is a bijection. Equation~(\ref{eq:wordinsertion})
follows from iteration of lemma~\ref{lemma:insertone},
where $L$ stems from $l$, $W$ from $w$, and $U$ from $u$.
\end{proof}

\begin{lemma}\label{lemma:sigma}
$$\sum_{w\in\SW_{i,j}}\sigma(w)=\left\{\begin{array}{lr}
0 & \mbox{if $i$ is even and $j$ is odd,} \\
(-1)^{\lfloor j/2\rfloor}
{{\lfloor i/2\rfloor}\choose{\lfloor j/2\rfloor}}
& \mbox{otherwise.}\end{array}\right.$$
\end{lemma}
\begin{proof}
By definition, we have $\sigma(w)=(-1)^L$, where $L=(w_1-1)+\cdots+(w_j-1)$.
Since $\sigma(w_1w_2\cdots w_j)\neq\sigma((w_1+1)w_2\cdots w_j)$
we only have to consider words in which $w_1+1=w_2$ and this value is
even. By iteration
of this argument we see that we only have to consider words
in which $w_{2k-1}+1=w_{2k}$ for $1\leq k\leq\lfloor j/2\rfloor$
and these values are even.
Every such pair gives an odd contribution to $L$.

If $j$ is odd, the last letter $w_j$ may be anywhere in the interval
$(w_{j-1},i]$. Since we have
$\sigma(w_1\cdots w_n)\neq\sigma(w_1\cdots(w_n+1))$
only words with $w_n=i$ odd remain. Then $w_n$ gives
an even contribution to $L$ so we can ignore it.

Thus, if $i$ is even and $j$ is odd the sum vanishes, otherwise
we can place the $\lfloor j/2\rfloor$ pairs
in $\lfloor i/2\rfloor$ positions, and we get
$(-1)^{\lfloor j/2\rfloor}
{{\lfloor i/2\rfloor}\choose{\lfloor j/2\rfloor}}$.
\end{proof}

\vspace{2mm}
\paragraph{\bf Remark.}
A referee has pointed out that, using $q$-binomial coefficients,
the sum in lemma~\ref{lemma:sigma} can be written
$$(-1)^{j\choose2}\left[
{\begin{array}{c}i \\ j\end{array}}\right]_{q=-1}.$$
This follows from the bijection between sorted words
$w_1w_2\cdots w_j\in\SW_{i,j}$ and weakly increasing sequences
$0\leq w_1-1\leq w_2-2\leq\cdots\leq w_j-j\leq i-j$, and from the
fact that $q$-binomial coefficients enumerate lattice paths
by area.
\vspace{2mm}

\begin{propo}\label{propo:horizontal}
Given an $n$-shape $\lambda$ whose strip consists of vertical dominoes, and
a nonnegative integer $k$,
let $H_\lambda$ be the set of $(n+k)$-shapes $\mu\supseteq\lambda$
with $\fb(\mu)=\fb(\lambda)$,
$\vs(\mu)=\vs(\lambda)$, and $\hs(\mu)=\lfloor k/2\rfloor$.
Then
$$\sum_{\mu\in H_\lambda}I_\mu
={{n/2+\lfloor k/2\rfloor}
\choose{\lfloor k/2\rfloor}}I_\lambda.$$
\end{propo}
\begin{proof}
Put $m=n+k$
and let $H_\lambda^\ast\supseteq H_\lambda$
be the set of $m$-shapes $\mu\supseteq\lambda$
with $\fb(\mu)=\fb(\lambda)$ and $\vs(\mu)=\vs(\lambda)$,
i.e.\ the set of $m$-shapes $\mu\supseteq\lambda$
with $v(\mu/\lambda)=0$.
By proposition~\ref{propo:colourbalance} all
$\mu\in H_\lambda^\ast\setminus H_\lambda$ have $I_\mu=0$.
Now we apply lemma~\ref{lemma:wordinsertion} to
$\T_{\lambda,k}$ and get
\begin{equation}\label{eq:Hlambda}
\sum_{\mu\in H_\lambda}I_{\mu}
=\sum_{\mbox{\scriptsize
$\begin{array}{c}T\in\T_{\lambda,k}\\ \sh(T\leftarrow w_{T,k})\in H_\lambda
\end{array}$}}\sigma(w_{T,k})\tau(\sh(T\leftarrow w_{T,k})/\lambda)\sgn(T).
\end{equation}
If $\sh(T\leftarrow w_{T,k})\in H_\lambda$ we have
$W=\lfloor k/2\rfloor$ (by the proof of
proposition~\ref{propo:colourbalance})
and $U$ is even in definition~\ref{defi:sigmatau},
which means that
$\tau(\sh(T\leftarrow w_{T,k})/\lambda)=(-1)^{\lfloor k/2\rfloor}$.
By first considering a summation of $\sigma(w_{T,k})\sgn(T)$ over
the whole set $H_\lambda^\ast$ and
then removing the contribution from $H_\lambda^\ast\setminus H_\lambda$,
we can write~(\ref{eq:Hlambda}) as
$$(-1)^{\lfloor k/2\rfloor}\left(\sum_{w\in\SW_{m,k}}\!\!\!\!\!\!\sigma(w)
\!\!\!\!\sum_{\mbox{\scriptsize
$\begin{array}{c}T\in\T_{\lambda,k}\\ w_{T,k}=w\end{array}$
}}\!\!\!\!\!\!\!\!\!\!\sgn(T)\ -
\sum_{\mu\in H_\lambda^\ast\setminus H_\lambda}
\!\!\!\!\sum_{\mbox{\scriptsize
$\begin{array}{c}T\in\T_{\lambda,k}\\
\sh(T\leftarrow w_{T,k})=\mu\end{array}$}}
\!\!\!\!\!\!\!\!\!\!\!\!\!\sigma(w_{T,k})\sgn(T)\right)$$
which equals
$$(-1)^{\lfloor k/2\rfloor}\left(\sum_{w\in\SW_{m,k}}\sigma(w)I_\lambda
-\sum_{\mu\in H_\lambda^\ast\setminus H_\lambda}
\frac{I_\mu}{\tau(\mu/\lambda)}\right)
=(-1)^{\lfloor k/2\rfloor}I_\lambda\sum_{w\in\SW_{m,k}}\sigma(w)$$
since $I_\mu=0$ for $\mu\in H_\lambda^\ast\setminus H_\lambda$.
By lemma~\ref{lemma:sigma},
$\sum_{w\in\SW_{m,k}}\sigma(w)=(-1)^{\lfloor k/2\rfloor}
{{n/2+\lfloor k/2\rfloor}\choose{\lfloor k/2\rfloor}}$
which gives the desired result.
\end{proof}

Proposition~\ref{propo:hooks} is now proved ``for free'':

\begin{proof}[\bf Proof of proposition~\ref{propo:hooks}]
If $h+v=\lfloor n/2\rfloor$,
applying proposition~\ref{propo:horizontal} to
$(1^{2v})$ and $k=n-2v$ yields the coefficient of $q^vx^h$:
$$\sum_{\mu\in H_{(1^{2v})}}I_\mu={{v+h}\choose{h}}I_{(1^{2v})}
={{v+h}\choose{h}}.$$
By proposition~\ref{propo:colourbalance},
the coefficient of $q^vx^h$ vanishes if $h+v\neq\lfloor n/2\rfloor$.
\end{proof}

For the proof of theorem~\ref{theo:jonas} we will need the following
observation.
\begin{lemma}\label{lemma:fourlingshape}
A non-empty fourling shape $D$ has zero sign-imbalance, $I_D=0$.
\end{lemma}
\begin{proof}
By lemma~\ref{lemma:chess} we only have to consider chess tableaux.
But there are no chess tableaux on a non-empty fourling shape since all outer
corners (squares
without neighbours below or to the right)
are black and the last number is even.
\end{proof}
\noindent
We will also need the following fundamental proposition.
\begin{propo}\label{propo:transpose}
For all shapes $\lambda$ we have
$$I_{\lambda'}=(-1)^{d(\lambda)}I_\lambda.$$
\end{propo}
\begin{proof}
Let $x=(r_x,c_x)$ and $y=(r_y,c_y)$ be two squares in $\lambda$ sorted so that
$r_x\leq r_y$.
After transposition $x$ becomes $(c_x,r_x)$ and $y$ becomes $(c_y,r_y)$
in $\lambda'$. The book permutation order between $x$ and $y$
is changed if and only if $r_x<r_y$ and $c_x>c_y$.
Thus $I_{\lambda'}=(-1)^p I_\lambda$, where
$p$ is the number of pairs $(x,y)$ of squares in $\lambda$ with
$x$ north-east of $y$.

Let $n$ be the number of squares in $\lambda$.
By proposition~\ref{propo:colourbalance} we can assume that
$\lambda$ has $\lfloor n/2\rfloor$ white squares.
Take any $n$-SYT $T$ on $\lambda$. For each number
$i$ in $T$, let $p_i$ be the number of north-east pairs
containing $i$ and a smaller number. It is easy to see that
if $i$ is in the square $(r,c)$ we have
$p_i=i-rc=(i+1)-(r+c+(r-1)(c-1))$, where
$r+c$ is odd if the square is white
and even if it is black, while $(r-1)(c-1)$ is odd if and only if
the square is the south-east corner of a fourling in the fourling body.
Thus,
$p=\sum_{i=1}^np_i\equiv\frac{n(n+3)}{2}+\lfloor n/2\rfloor+d(\lambda)\pmod2$,
since there are $\lfloor n/2\rfloor$ white squares in $\lambda$.
But $\frac{n(n+3)}{2}+\lfloor n/2\rfloor=\lfloor n(n+4)/2\rfloor$
is always even, so $p\equiv d(\lambda)\pmod2$.
\end{proof}

Finally we have all the tools we need.
\begin{proof}[\bf Proof of theorem~\ref{theo:jonas}]
By proposition~\ref{propo:colourbalance},
we can assume that $h+v=\lfloor s/2\rfloor$.
Let $V$ be the set of shapes with fourling body $D$,
$2v$ squares in the strip, and $v$ vertical strip dominoes.
First we will show that $\sum_{\lambda\in V}I_\lambda=0$.
Let $V'=\{\lambda':\lambda\in V\}$.
By proposition~\ref{propo:transpose},
$\sum_{\lambda\in V}I_\lambda=(-1)^{d(D)}\sum_{\lambda\in V'}I_\lambda$,
so it suffices to show that the latter sum vanishes.
Applying proposition~\ref{propo:horizontal} to $D'$ and
$k=2v$ yields
$$\sum_{\lambda\in V'}I_\lambda=\sum_{\lambda\in H_{D'}}I_\lambda
={{2d(D)+v}\choose{v}}I_{D'}=0$$
by lemma~\ref{lemma:fourlingshape}.
Finally, we apply proposition~\ref{propo:horizontal} to every
$\lambda\in V$ and $k=s-2v$, and get
$$\sum_{\lambda\in V}\sum_{\mu\in H_\lambda}I_\mu
={{2d(D)+v+h}\choose{h}}\sum_{\lambda\in V}I_\lambda=0.$$
\end{proof}

\section{Possible generalizations}
\noindent
The concept of sign-imbalance generalizes naturally to general finite
posets. Note that a SYT is a linear extension of the
partial order on the squares implied by coordinate pairs.

Let $P$ be an $n$-element poset and let
$\omega:P\rightarrow[n]=\{1,2,\ldots,n\}$ be a bijection called
the {\em labelling} of $P$.
A {\em linear extension} of $P$ is an
order preserving bijection $f:P\rightarrow[n]$.
If we regard $f$ as a permutation $\pi_f$ of $[n]$
given by $\pi_f(i)=\omega(f^{-1}(i))$ we can talk about the
sign of $f$. The {\em sign-imbalance} of $P$ is the sum of the
signs of all linear extensions of $P$. If the sign-imbalance of $P$
is zero we say that $P$ is {\em sign-balanced}.

Note that the sign of a linear extension depends on the labelling $\omega$.
However, this dependence is not essential since changing the labelling
of $P$ simply multiplies $\pi_f$ by a fixed permutation.
For instance,
the sign-imbalance of $P$ is defined up to a sign without specifying $\omega$,
and the notion of sign-balance is completely independent of the labelling.

There has been some work (see \cite{stanley})
considering sign-imbalances of general posets
and identifying the sign-balanced ones.
Unfortunately, the approach taken in this paper does not seem
applicable to this more general question.

If we specialise to partition shapes, however, we
hope that our Robinson-Schensted technique will be useful
in future research. Some things to do:
\begin{itemize}
\item Characterise the sign-balanced partition shapes. There are some
theorems on sign-balanced posets (see~\cite{stanley}); a
complete characterization in the special case of partition shapes
may shed some light on this more general question.
\item Find the ``best'' version of theorem~\ref{theo:jonas},
i.e.\ find the smallest classes of $n$-shapes whose imbalance
sum vanishes. This is a generalization of the above and, as
figure~\ref{fig:jonasex} shows, there is still work to do.
\item Find a nice formula for $I_\lambda$, maybe in
the same spirit as the hook
length formula. This may very well be impossible, as Stanley points
out~\cite[page 14]{stanley}.
\item Study the imbalance of skew partitions. This is an
interesting issue since most structural properties of partitions
generalize to skew partitions, including the RS-correspondence (see
e.g.~\cite{saganstanley}).
\end{itemize}

\end{document}